\documentclass{svproc}
\usepackage{url}

\usepackage{amsmath,bm}
\usepackage{amssymb}
\usepackage{amsfonts}
\usepackage{subcaption}
\usepackage{tikz}
\usepackage[europeanvoltages,
  europeancurrents,
  europeanresistors,
  europeaninductors,
  smartlabels]{circuitikz}
\usepackage{pgfplots}
\usepackage{pgfplotstable}
\usepackage{tudacolors}
\usepackage{siunitx}

\usepackage{todonotes}
\usepackage{xcolor}
\definecolor{cb-black}      {RGB}{  0,   0,   0}
\definecolor{cb-blue-green} {RGB}{  0,  073,  073}
\definecolor{cb-green-sea}  {RGB}{  0, 146, 146}
\definecolor{cb-rose}       {RGB}{255, 109, 182}
\definecolor{cb-salmon-pink}{RGB}{255, 182, 119}
\definecolor{cb-purple}     {RGB}{ 73,   0, 146}
\definecolor{cb-blue}       {RGB}{ 0, 109, 219}
\definecolor{cb-lilac}      {RGB}{182, 109, 255}
\definecolor{cb-blue-sky}   {RGB}{109, 182, 255}
\definecolor{cb-blue-light} {RGB}{182, 219, 255}
\definecolor{cb-burgundy}   {RGB}{146,   0,   0}
\definecolor{cb-brown}      {RGB}{146,  73,   0}
\definecolor{cb-clay}       {RGB}{219, 209,   0}
\definecolor{cb-green-lime} {RGB}{ 36, 255,  36}
\definecolor{cb-yellow}     {RGB}{255, 255, 109}
\newcommand{\xk}{x_k}
\newcommand{\zk}{z_k}
\newcommand{\vk}{v_k}
\newcommand{\ik}{i_k}
\newcommand{\vold}{v_{k-1}}
\newcommand{\iold}{i_{k-1}}

\newcommand{\dkx}{\Delta_k x}

\newcommand{\dkv}{\Delta_k v}
\newcommand{\dki}{\Delta_k i}

\newcommand{\diold}{\Delta_{k-1}i}

\begin{document}
\mainmatter              
\title{Waveform Relaxation for Field/Circuit Coupled DAEs with Generalized Capacitances}
\titlerunning{Waveform Relaxation for Circuits with Generalized Capacitances}  
%
\author{Idoia Cortes Garcia\inst{1} \and Jonas Pade\inst{2}}
\authorrunning{I. Cortes Garcia and J. Pade} 
%
%
\institute{Department of Mechanical Engineering-Dynamics and Control, Eindhoven
University of Technology, Eindhoven, The Netherlands, \email{i.cortes.garcia@tue.nl}
\and
Humboldt-Universit\"at zu Berlin, Department of Mathematics,  10099 Berlin, Germany, \email{jonas.pade@hu-berlin.de}}

\maketitle              

\begin{abstract}
Field/circuit coupling is a common approach when a lumped representation of a certain electrotechnical device is not accurate enough. To exploit existing code and underlying properties of the coupled systems, cosimulation techniques such as waveform relaxation can be used. The coupled system is of differential-algebraic type, which can potentially lead to divergence. This paper presents a novel, sufficient topological convergence criterion for field/circuit coupled systems of higher index containing a generalized capacitance. Hereby, the criterion holds for a full range of field systems whose structure can be classified as a generalized capacitance.
Finally, the theoretical results are supported by numerical simulations.
\end{abstract}

\section{Introduction}
Lumped circuit are ubiquitous in electrical engineering and often used to describe the electromagnetic behaviour of engineering systems. However, they neglect spatially distributed electromagnetic phenomena, which, for some applications, is too inaccurate.
In these cases, field/circuit coupling can be performed \cite{Tsukerman_1993aa,Bartel}, where fully spatially distributed electromagnetic field elements are coupled to a surrounding lumped circuit. Recent work has shown that these field models can be understood as a generalization of classical circuit elements \cite{CortesGarcia_2020}.

Cosimulation methods are well suited to solve field/circuit coupled systems since they allow solving the field and circuit parts separately and thus the usage of dedicated solvers and time scales for the different subsystems \cite{Bartel}.
Originally introduced in \cite{Lelarasmee}, a well-established family of iterative cosimulation methods is waveform relaxation (WR), also known as dynamic iteration. {Its application and convergence behaviour for differential-algebraic equations (DAEs) is still a topic of research \cite{Bartel_2023aa}.} While WR is convergent for coupled ODEs on finite time intervals \cite{Burrage},
it can diverge for coupled DAEs unless an additional contractivity condition is satisfied \cite{Arnold}.
Computing the contractivity condition, however, is usually costly. 
In this work, we extend the easy-to-check topological WR convergence results for higher-index DAEs obtained in \cite{Pade_diss} to
coupled field/circuit systems whose field is a generalized capacitance \cite{CortesGarcia_2020}. 
With this, the theoretical results hold for all distributed models with a similar structure, that can be classified as a generalized capacitance.
If the convergence criterion is not satisfied, we have to expect WR to diverge or to converge comparatively slowly. 

The paper is structured as follows. Section~\ref{sec:pade_coupl} introduces the definiton of generalized capacitances and the field/circuit coupled system. The WR scheme and main convergence results are presented in Section~\ref{sec:pade_conv}. Finally, Section~\ref{sec:padeexample} and~\ref{sec:padeconclusions} finalize the paper with numerical simulation examples and conclusions.

 \section{Field/Circuit Coupled Model}
 \label{sec:pade_coupl}
Here, we introduce the coupled field/circuit system whose field element can be classified as a so-called \emph{generalized capacitance} as in \cite{CortesGarcia_2020}.
\subsection{Generalized Capacitance}
A classic, lumped capacitance is described by the current-to-voltage relation
\begin{align}
\label{eq:lumped_cap}i=C(v)v'.
\end{align}
The state-dependent capacitance $C(\cdot)$ is required to be nonsingular for all states $v$. This ensures that \eqref{eq:lumped_cap} can be solved for $v'$.
That way, for a given current excitation $i$, the voltage $v$ is described as an ODE $v'=C^{-1}(v)i$ which depends on the input $i$.
This well-posed ODE nature of the capacitance model w.r.t. unknown $v$ and input $i$ should be reflected in a generalized definition of a capacitance.
\begin{definition}\label{def:padedef}
A {\textbf {generalized capacitance}} is described by
\begin{align}
     \label{eq:genel}f\left(\frac{d}{dt}m(x,i,v,t),x,i,v,t\right)=0,
\end{align}
where $f$ is continuously differentiable and after at most one differentiation, we can reformulate  \eqref{eq:genel}  with locally Lipschitz continuous functions $\chi$ and $g$ and a continuously differentiable function $R_\chi$ as
\begin{align}
 \label{eq:gencapx} x' &= \chi(x,i,v,t)+R_\chi(i)i'\\
 \label{eq:gencapv} v' &= g(x,i,v,t).
\end{align}
\end{definition}
In the case of electromagnetic field models, the unknown $x$ usually stems from spatial discretization of the field.
\begin{remark}
The above generalized definition can be understood as an index-1 condition for the unknowns $(x,v)$ with input $i$. Additional robustness is required for $v$, which manifests itself in the absence of a derivative $i'$ in  \eqref{eq:gencapv}.  We emphasize that we obtain an index 1 DAE in unknowns $(x,v)$ for input $i$, and higher index in $(x,i)$ for input $v$.
In a coupled model, the network topology determines whether or not this leads to a higher index network DAE.
\end{remark}
\begin{remark}
In \cite{CortesGarcia_2020}, "strongly capacitance-like" elements are defined, satisfying certain strong monotonicity conditions. Physically, these are related with passivity properties of the elements. From a mathematical point of view, they make sure that the index results from \cite{topal} hold for generalized elements, too (c.f. \cite{CortesGarcia_2020}). Since these conditions and the implied index results from \cite{topal} are not necessary for our convergence analysis, and quite technical, we do not include them here.
\end{remark}
 
\subsection{Circuit and Coupled Model}
Generalized capacitances are designed to be embedded in a surrounding circuit. In the following, its coupling is presented.

The circuit subsystem is modeled by means of the \textit{modified nodal analysis} (MNA) \cite{mna} and it consists of quasilinear inductances and capacitances, fully nonlinear diodes and resistances, and independent voltage and current sources. In compact notation and assuming the circuit is coupled to a field element via the voltage $v$ and the current $i$, it reads 
 \begin{align}\label{eq:circ}
E(z) z'+f(t,z)=Pi,\qquad P^\top z = v.
\end{align}
Here, the unknown $z$ consists of node potentials and currents through inductors and voltage sources, $E$ and $f$ contain the system equations arising from MNA and 
$P^\top=(A_{\mathrm{f}}^\top\ 0)$, where $A_{\mathrm{f}}$ is the incidence matrix of the field element. We assume certain common local Lipschitz and passivity properties, cf. \cite{topal} for details.


The MNA circuit system coupled to a generalized capacitance field element with coupling variables $i$ and $v$ reads
\begin{align}
\label{eq:coupfield} f\left(\frac{d}{dt}m(x,i,v,t),x,i,v,t\right)&=0,\\
 \label{eq:coupcirc}
E(z) z'+f(t,z)=Pi,\qquad P^\top z &= v.
\end{align}

\section{Waveform Relaxation Scheme for Coupled Model}
\label{sec:pade_conv}
For cosimulation we consider the Jacobi WR scheme. In our coupled model, it reads
\begin{align}
\label{eq:wrfield} f\left(\frac{d}{dt}m(\xk,\iold,\vk,t),\xk,\iold,\vk,t\right)&=0,\\
\label{eq:wrcirc} E(\zk) \zk'+f(t,\zk)=Pi_k,\qquad P^\top z_k&=\vold,
\end{align}
where $k$ is the iteration parameter.
Note that we chose a voltage driven circuit and a current driven field here, which is reflected by the voltage $v$ and the current $i$ entering the circuit  \eqref{eq:wrcirc} and the field  \eqref{eq:wrfield} as
previous $(k-1)$-th iterates, hence inputs, respectively.

Note that it is possible to choose the  generalized capacitance as voltage driven, but not advisable  from the WR point of view. Notably, this would
render an index-1 type estimate as in  \eqref{eq:padep} impossible, which is crucial for convergence Theorem \ref{theo:convergence}.

\paragraph{Consistent initial values}
Initial values for DAEs are only partly (mathematically) free, since some are fixed by algebraic constraints.
We bypass the practical difficulty of finding consistent initial values (those which satisfy the algebraic constraints of DAEs) as follows: For a general DAE $f(x',x,t)=0$, we define $x_0^{dyn}=Bx(t_0)$ as the set of dynamic, hence freely selectable initial values, where $B$ is an appropriate matrix (which is generally non-trivial to find and can be nonlinear).
We refer to \cite{Estevez-Schwarz_2000ab} for details on consistent initalization for DAEs.

\subsection{A Convergence Theorem}
We call two nodes CV-connected, if there exists a path of (no elements other than) capacitances and voltage sources between them. In other words, there exists no LIR-cutset between such nodes.
\begin{theorem}\label{theo:convergence}
Let the coupling nodes of the field/circuit system \eqref{eq:coupfield}-\eqref{eq:coupcirc} be not CV-connected. Then, the WR scheme \eqref{eq:wrfield}-\eqref{eq:wrcirc} is convergent. 

 More precisely:
 For any given consistent initial values on a sufficiently small time interval $\mathcal I=[t_0,T]$, the sequence of solutions of the Jacobi WR scheme \eqref{eq:wrfield}-\eqref{eq:wrcirc} 
 converges for $k\to\infty$ in $(C(\mathcal I),\|\cdot\|_\infty)$ to the solution of the corresponding field/circuit coupled model \eqref{eq:coupfield}-\eqref{eq:coupcirc}. The rate of convergence is $\sqrt{Hc}$ for an appropriate constant $c$ and $H:=T-t_0$.
\end{theorem}
\begin{remark}
Due to the beneficial effect of a small time interval $H$ on WR convergence, so-called \emph{windowing} is common for WR methods. This involves the partition of the time interval into smaller subintervals and the successive execution of WR on each subinterval separately.
\end{remark}
\begin{remark}
 An analogous result holds for Gauss-Seidel WR, where the rate of convergence is $Hc$. The faster convergence is a tradeoff against parallelization, which is in turn only possible for Jacobi WR. 
\end{remark}
\subsection{Proof of the Convergence Theorem}
For this subsection, we introduce the notation
\begin{align*}
    \Delta_k \alpha:=\alpha_k-\alpha,\qquad \Delta_k \alpha_0:=\alpha_k(t_0)-\alpha(t_0),\qquad \alpha_0:=\alpha(t_0)
\end{align*}
for any variable $\alpha$. Furthermore, we denote by $|\cdot|$ the maximum norm in $\mathbb R^n$ and by $\|\cdot\|:=\|\cdot\|_\infty$ the uniform norm for continuous functions on $\mathcal I$.

\begin{lemma}\label{lem:cvcirc} 
Consider the circuit \eqref{eq:circ} with a continuous perturbation $\delta$ of $v$:
\begin{align}\label{eq:padepertcirc}
E(z_\delta) z_\delta'+f(t,z_\delta)=Pi_\delta,\qquad P^\top z_\delta = v+\delta.
\end{align}
 If the coupling nodes are not CV-connected, then there exists a constant $\nu$ such that for the difference of solutions $\Delta(z,i):=(z_\delta-z,i_\delta-i)$  of the perturbed and the unperturbed circuit, it holds for all $t\in \mathcal I$
\begin{align}\label{eq:padep}
 |\Delta(z,i)(t)|\leq \nu|\delta(t)|
\end{align}
\end{lemma}
For a proof of the Lemma, we refer to \cite{Pade_diss}, pp.55, Theorem 3.5.1. 
\begin{remark}
The well-known topological index analysis introduced in \cite{topal} identifies CV-loops (and LI-cuts) in circuits as the cause for resulting index-2 DAEs. In this context, the lemma's statement can be reformulated as follows: a local perturbation of a voltage source has an index-1 type impact (as described by \eqref{eq:padep}) on the state of the circuit. Note that this is independent of the (global) DAE index of  \eqref{eq:padepertcirc}. This is the crucial idea for the proof, which employs the well-known topological index concept \emph{only locally on the coupling variables}. They
are decisive for the convergence of WR. In practice, this means that also systems with index 2 or higher can lead to a convergent co-simulation, as long as the coupling variables have a lower-index relation.
\end{remark}
\paragraph{Proof of Theorem}
For the proof we require the estimates:
\begin{enumerate}
   \item $\|\dki\|\leq \eta\|\Delta_{k-1}v\|$
  \item $\|\dkx\|\leq H\lambda_1\|\dkv\|+\lambda_2\|\diold\|$
   \item $\|\dkv\|\leq H\mu\|\diold\|$
    \end{enumerate}
Plugging in 1. into 3., we then obtain the recursion 
\begin{align}
\label{eq:padea}\|\dkv\|\leq Hc\|\Delta_{k-2}v\|, \qquad c=\mu\eta
\end{align}
which converges to $0$ with rate $\sqrt{Hc}$ for $k\to \infty$ if $Hc<1$.
Due to 1. and 2., this implies that $\Delta_k i$ and 
$\Delta_k x$ converge with the same rate. Considering Lemma \ref{lem:cvcirc}, Estimate \eqref{eq:padea} furthermore implies convergence of $\Delta_kz$, again with rate $\sqrt{Hc}$. In the following we derive estimates 1-3. 

1.
We notice that solving the circuit subsystem \eqref{eq:wrcirc} of the WR scheme is equivalent to solving the circuit with a voltage perturbation as in \eqref{eq:padepertcirc} with $\delta=v_{k-1}-v$.
The desired estimate is then immediately implied by Lemma \ref{lem:cvcirc}. 

 2. For the field subsystem, integration of 
 \eqref{eq:gencapx}-\eqref{eq:gencapv} yields
 \begin{align}\label{eq:intform}
  x(t)=x_0+\int_{t_0}^t\chi(x,i,v,\tau)+R_\chi(i)i'd\tau,\qquad
  v(t)=v_0+\int_{t_0}^t g(x,i,v,\tau)d\tau,
 \end{align}
 where we drop the argument $\tau$ of the functions $x,i,v$ within the integral for better readability.

 Writing the analogous integral formulation for the corresponding equations of the WR scheme and with the difference $x_{k}-x=\Delta_kx$, basic calculus yields
 \begin{align*}
  \left|\int_{t_0}^t\chi(\xk,\iold,\vk,\tau)-\chi(x,i,v,\tau)d\tau\right|
  \leq L_\chi\left(\int_{t_0}^t|\dkx|d\tau+|\diold|+|\dkv|d\tau\right)
 \end{align*}
and, denoting the integral of $R_\chi(\cdot)$ by $r_\chi(\cdot)$,
 \begin{align*}
  \left|\int_{t_0}^tR_\chi(\ik)(\ik')-R_\chi(i)i'd\tau\right|&\leq \left|r_\chi(\ik(t))-r_\chi(i(t))\right|+\left|r_\chi(i_{k_0})-r_\chi(i_0)\right|\\
  &\leq L_{r_\chi}\left(\left|\dki(t)|+|\dki_0\right|\right),
 \end{align*}
where $L_*$ are the respective Lipschitz constants whose existence on bounded sets is guaranteed by definition of a generalized capacitance.

 Starting from the first equation in \eqref{eq:intform} and considering the estimates thereafter, the lemma of Gronwall yields the existence of a constant $\tilde c$ such that for all $t\in\mathcal I$
\begin{align}
 |\dkx(t)|&\leq \tilde c\left(\int_{t_0}^t|\diold|+|\dkv|d\tau\right)+L_{r_\chi}\left(|\diold(t)|+|\diold_0|\right)\nonumber\\
 &\leq H\lambda_1\|\dkv\|+\lambda_2\|\diold\|.\label{eq:xest}
\end{align}
Note that the constant from the Gronwall lemma has the form $c(H)=e^{\beta H}$ for a constant $\beta$. However, since we assume $H$ sufficiently small hence bounded, there exists $\tilde c\geq  c(H)$ independent of $H$ which satisfies \eqref{eq:xest}.

3. The second equation in \eqref{eq:intform} gives for all $t\in\mathcal I$
\begin{align*}
 |\dkv(t)|\leq \int^t_{t_0}|\dkx|+|\diold|d\tau
\end{align*}
hence, after plugging in \eqref{eq:xest} to eliminate $\dkx$, the lemma of Gronwall yields the existence of a $\hat c$ such that
\begin{align*}
 \|\dkv\|\leq H^2\hat c\|\diold\|+H\hat c(1+\lambda_2)\|\diold\|\leq H\mu\|\diold\|.
\end{align*}
\qed
\section{Numerical Examples}
\label{sec:padeexample}
\begin{figure}
\centering
\begin{subfigure}{0.45\textwidth}
\centering
\begin{tikzpicture}
\def \len {2.2}
\draw (0,0) to[V=$v_{\mathrm{s}}$, *-*] (\len,0)
    to [C=$C$, -*] (\len,\len)
    to [R=$R$] (0,\len)
    to (0,0) -- (0,-0.15) node[ground]{}; 
\draw (\len,0) to[I=$i_{\mathrm{s}}$, -*] (2*\len, 0)
    to[L=$L_2$,-*] (2*\len, \len)
    to[L=$L_1$] (\len,\len);
\draw (0,\len) -- (0,\len+0.5*\len)
    -- (2*\len,\len+0.5*\len)
    -- (2*\len,\len);
\node [draw, thick, fill=white, shape=rectangle, minimum width=1.5cm, minimum height=0.5cm, anchor=center] at (\len,\len+0.5*\len) {\includegraphics[width=1cm]{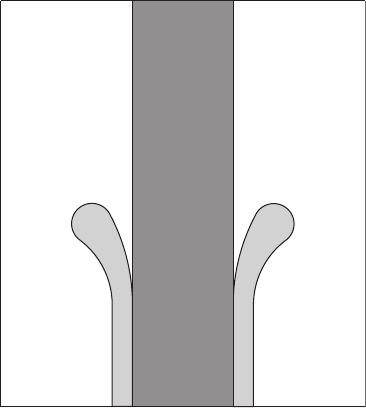}};
\end{tikzpicture}
\end{subfigure}
\begin{subfigure}{0.45\textwidth}
\centering
\begin{tikzpicture}
\def \len {2.2}
\draw(0,0) to[V=$v_{\mathrm{s}}$, *-*] (\len,0)
    to [C=$C$, -*] (\len,\len)
    to [R=$R$] (0,\len)
    to (0,0) -- (0,-0.15) node[ground]{}; 
\draw (\len,0) to[I=$i_{\mathrm{s}}$, -*] (2*\len, 0)
    to[L=$L_2$,-*] (2*\len, \len)
    to[L=$L_1$] (\len,\len);
\draw (0,\len) -- (0,\len+0.5*\len)
    -- (\len,\len+0.5*\len)
    -- (\len,\len);
\node [draw, thick, fill=white, shape=rectangle, minimum width=1.5cm, minimum height=0.5cm, anchor=center] at (0.5*\len,\len+0.5*\len) {\includegraphics[width=1cm]{figures/cable_terminator_drawing.pdf}};
\end{tikzpicture}
\end{subfigure}
\caption{Convergent (first) and divergent (second) circuit couplings. The circuits and toy values for the parameters ($R=1\unit{\ohm}$, $L_1 =5\unit{\henry}$, $L_2=2\unit{\henry}$, $C=1\unit{\farad}$)  are based on examples
of \cite{Pade_2021aa}.}
\label{fig:circuits}
\end{figure}

In the following we illustrate the convergence results in Section~\ref{sec:pade_conv} with a toy example. For the field subsystem we choose the electroquasistatic description of the cable termination model in \cite{CortesGarcia_2019aa}, which is a capacitance-like element as proven in \cite{CortesGarcia_2020}. We choose the surrounding circuits of \cite{Pade_2021aa}, with an additional inductance to increase the system's index to $2$. The circuits and its coupling with the field model are sketched in Figure~\ref{fig:circuits}. The waveform relaxation cosimulation is done with a Gauss-Seidel WR scheme. 

\begin{figure}
\begin{tikzpicture}
\begin{axis}[width = 0.47\textwidth,
        xlabel=$t$,
        ylabel=Voltage $v_{\mathrm{f}}$,
         y label style={at={(0.1,0.5)}},
        legend style={at={(1.07,0.95)},anchor=west}]
\addplot[line width = 1pt, cb-rose] table [x=t, y=v, col sep=comma] {simulation_data/circb_monolithic.csv};
\addplot[line width = 1pt, dashed, black] table [x=t, y=x1, col sep=comma] {simulation_data/circb_conv_solution_field.csv};
\addplot[line width = 1pt, dashed, cb-blue-green] table [x=t, y=x2, col sep=comma] {simulation_data/circb_conv_solution_field.csv};
\addplot[line width = 1pt, dashed, cb-burgundy] table [x=t, y=x3, col sep=comma] {simulation_data/circb_conv_solution_field.csv};
\legend{$v_{\mathrm{mono}}$,$v_{\mathrm{f}}^1$,$v_{\mathrm{f}}^2$,$v_{\mathrm{f}}^3$}
\end{axis}
\hspace{0.57\textwidth}
\begin{axis}[width = 0.47\textwidth,
        xlabel=$t$,
        ylabel=Voltage $v_{\mathrm{f}}$,
        legend style={at={(0.7,0.25)},anchor=west},
       y label style={at={(0.15,0.5)}},
        domain = -1.5:1.5]
\addplot[line width = 0.5pt, dashed, cb-burgundy, forget plot] table [x=t, y=x3, col sep=comma] {simulation_data/circa_div_solution_field.csv};
\addplot[line width = 0.5pt, dashed, cb-blue-green, forget plot] table [x=t, y=x2, col sep=comma] {simulation_data/circa_div_solution_field.csv};
\addplot[line width = 0.5pt, dashed, black] table [x=t, y=x1, col sep=comma] {simulation_data/circa_div_solution_field.csv};
\addplot[line width = 1pt, cb-rose] table [x=t, y=v, col sep=comma] {simulation_data/circa_monolithic.csv};
\end{axis}
\end{tikzpicture}
\caption{Evolution of voltage across generalized capacitance of monolithic  $v_{\mathrm{mono}}$ and co-simulated $v_{\mathrm{f}}^k$, $k=1,\ldots,8$ solutions for convergent (left) and divergent (right) circuits in Figure~\ref{fig:circuits}.}\label{fig:wr}
\end{figure}
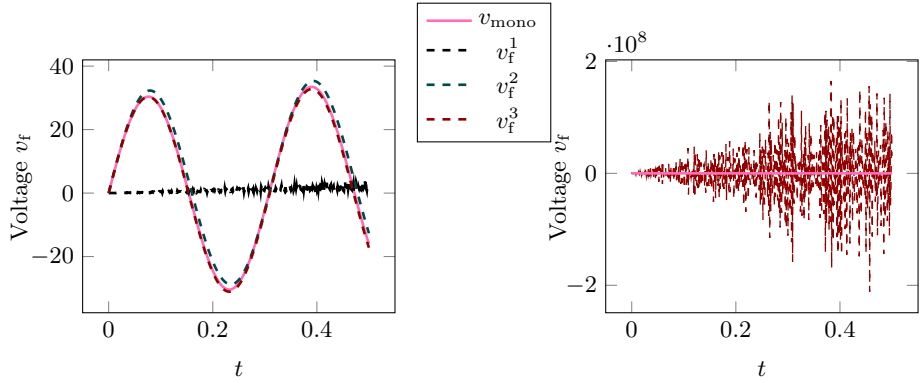
\begin{figure}
\centering
\begin{tikzpicture}
    \begin{axis}[width = 0.5\textwidth,
                xmode=log,
                ymode=log,
                xlabel = $H$,
                ylabel = Error,
                legend style={at={(0.7,0.25)},anchor=west}
                ]
        \addplot[line width = 1pt, black, mark=*] table [x=H, y=errV, col sep=comma] {simulation_data/convergence_cirb.csv};
        \addplot[line width = 1pt, gray] table [x=H, y expr={300*(\thisrowno{0})}, col sep=comma] {simulation_data/convergence_cirb.csv};
        \legend{$\mathrm{err}_{4,\mathrm{v}}$,$\mathcal{O}(H)$}
    \end{axis}
\end{tikzpicture}
\caption{Error $\mathrm{err}_{\mathrm{v}} = \max_t\Delta_1v(t)$ for different windows of size $H$.}\label{fig:convergence}
\end{figure}
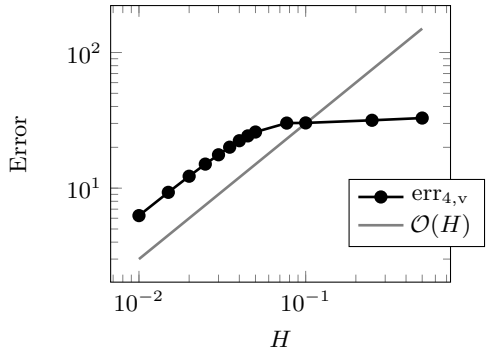

Figure \ref{fig:wr} shows the voltage $v_{\mathrm{f}}^k$ across the field element on the $k$-th WR iteration for both circuit coupling cases of Figure~\ref{fig:circuits}. As proven by theory, the WR scheme is convergent for the first coupled circuit, and divergent for the second one. As expected, this is the opposite convergence behaviour as the results in \cite{Pade_2021aa}, where an inductance-like element is coupled to analogous circuits.
Furthemore, Figure~\ref{fig:convergence} shows the convergence order for the convergent case for fixed $k=1$ and different window sizes $H$. As expected, our Gauss-Seidel scheme converges with order $H$ for small enough window sizes (after some preasymptotic behaviour).

\section{Conclusions}
\label{sec:padeconclusions}
We have presented a topological convergence criterion for WR schemes of circuits coupled to electromagnetic fields, where the field can be classified as a generalized capacitance.  The criterion that the coupling nodes must not be CV-connected holds for all such coupled field/circuit systems. Simulations of two circuit examples coupled to an electroquasistatic cable termination model confirm the theoretical results. We believe that similar convergence conditions can be derived for generalized inductances and resistances to cover a wide range of distributed systems. This, however, is left for future research.

\paragraph{\bf Acknowledgement} This work was funded by the DFG, Germany, TRR 154, Mathematical Modelling, Simulation and Optimization Using the Example of Gas Networks, project number 239904186, project C02.

\par\vfill
\end{document}